\documentclass[11pt]{amsart}
\usepackage{amssymb}
\usepackage[english]{babel}
\usepackage{mathrsfs}

\newtheorem{theorem}{Theorem}[section]
\newtheorem{lemma}[theorem]{Lemma}
\newtheorem{proposition}[theorem]{Proposition}
\newtheorem{corollary}[theorem]{Corollary}

\theoremstyle{definition}
\newtheorem{definition}[theorem]{Definition}
\newtheorem{example}[theorem]{Example}
\newtheorem{problem}[theorem]{Problem}
\newtheorem{remark}[theorem]{Remark}
\newtheorem{claim}{Claim}
\newtheorem{fact}[theorem]{Fact}

\def\N{\mathbb{N}}
\def\P{\mathbb{P}}
\def\Prm{\P}
\def\Z{\mathbb{Z}}

\def\cont{\mathfrak{c}}
\def\cd{\text{$|\cdot|$}-divisible}

\begin{document}
	
	\title[Algebraic structure of pseudocompact abelian groups]{A complete description of the algebraic structure of pseudocompact
	abelian groups}
	
\author[D. Dikranjan]{Dikran Dikranjan}
\address{Department of Mathmatical, Computer and Physical Sciences,
University of Udine\\
Via delle Scienze 106\\
33100 Udine, 
Italy}	
\email{dikran.dikranjan@uniud.it}

	\author[A. Dorantes-Aldama]{Alejandro Dorantes-Aldama}
	\address{Universidad Panamericana Campus M\'exico, Augusto Rodin 498, Ciudad de M\'exico, 03920, M\'exico}
	\email{adorantes@up.edu.mx}

	\author[D. Shakhmatov]{Dmitri Shakhmatov}
 	\address{
Graduate School of Science and Engineering\\
		Ehime University, Matsuyama 790-8577, Japan}
	\email{dmitri.shakhmatov@ehime-u.ac.jp}
	\thanks{The second and the third listed author were partially supported by the Grant-in-Aid for Scientific Research~(C) No.~20K03615
	 by the Japan Society for the Promotion of Science (JSPS)}

\dedicatory{In memory of W. Wistar Comfort (1933-2016)}
	
	\begin{abstract} 
We solve a major open problem 
from \cite{DiS}
by describing abelian groups which admit a pseudocompact Hausdorff group topology. 
First, we prove that an abelian group 
$G$ 
can be equipped with a pseudocompact Hausdorff group topology iff it can be decomposed into a direct 
sum $G=B\oplus H$ of a bounded torsion group $B$ and  
a
group $H$ satisfying $|nH|=|H|$ for all positive integers $n$, such that both $B$ and $H$ admit a pseudocompact Hausdorff group topology.
Second, for an unbounded abelian group $G$, we denote by $m$ the minimal positive integer $n$ 
such that $|nG|=|kG|$ for all integers $k\ge n$.
Then $G$ admits a pseudocompact group topology iff
both $G[m]=\{g\in G: mg=0\}$ and $mG\cong G/G[m]$ admit a pseudocompact Hausdorff group topology.
	\end{abstract}
	
	\maketitle
	
{\em All topological groups and group topologies are assumed to be Hausdorff\/}.

	We use $\N$ to denote the set of {\em positive\/} integers and $\Prm$ to denote the set of prime numbers.

Let $G$ be an abelian group and $n\in\N$.
Note
that
$nG=\{ng:g\in G\}$
is a subgroup of $G$ and the map
$g\mapsto ng$ ($g\in G$) is a homomorphism of $G$ onto $nG$.
Similarly, 
$G[n]=\{g\in G: ng=0\}$
is a subgroup of $G$ that 
is closed in every 
(Hausdorff) group topology on $G$.
The subgroup $nG$ of $G$ is algebraically isomorphic to the quotient $G/G[n]$ of $G$.

An abelian group $G$ is: 
\emph{torsion} if $G=\bigcup\{G[n]:n\in\N\}$, and
\emph{bounded (torsion)} if $G=G[n]$, or equivalenly $nG=\{0\}$, for some $n\in\N$.
The minimal 
$n\in\N$ such that $nG=\{0\}$
is the {\em order\/} of a bounded group $G$ and shall be denoted by $n_G$. 
For
$p\in\Prm$, 
$G$ is called a {\em $p$-group\/} provided that
$G = \bigcup \{G[p^n] : n \in\N\}$.

Notation $G\cong H$ means that abelian groups are isomorphic.

For a cardinal $\tau$ we denote by $G^{(\tau)}$ the direct sum of $\tau$ copies
of the group $G$.

For a topological space $X$ we use $w(X)$ to denote the {\em weight\/} of $X$ defined as
the minimum cardinality of a base for $X$. Note that a finite Hausdorff space $X$ is discrete, so satisfies $w(X)=|X|$.

\section{Introduction}

A Tychonoff space is \emph{pseudocompact} if every real-valued continuous function on it is bounded.
Compact spaces are pseudocompact.

A topological group $G$ is \emph{precompact} or \emph{totally bounded} provided that 
$G$ is a subgroup of some compact group.
For every topological group, the following holds:
\begin{equation}
\label{compact:implies:pseudocompact}
\text{compact}\to \text{pseudocompact}\to \text{precompact}.
\end{equation}

Both the class of compact groups and the class of precompact groups are closed under arbitrary direct products, thanks to the Tychonoff product theorem.
The product of even two pseudocompact spaces need not be pseudocompact.
Comfort and Ross  showed 
that pseudocompactness becomes productive for topological groups:
\begin{theorem}
\label{Comfort-Ross-theorem}
\cite{CRoss}
The class of pseudocompact groups is closed under arbitrary direct products.
\end{theorem}

Every abelian group can be equipped with a (Hausdorff) precompact group topology.
Equivalently, every discrete abelian group is maximally almost periodic (MAP).
The algebraic structure of compact abelian groups is well-known \cite{Hulanicki1, Hulanicki2, Harrison}.
In this paper, we shall 
completely resolve the 
following long-standing problem:

\begin{problem}\cite[Problem 0.2]{DiS}
\label{main:problem}
Describe abelian groups which can be equipped with a (Hausdorff) pseudocompact group topology.
\end{problem}
A lot of activity 
attempting to resolve this problem happened in the 1990s \cite{vD, DS-Nish, DS-York, DS-TopProc, CR1, CR2, DiS}.
Below we summarize the most important known results about the problem. 
\begin{definition}\cite[Definition 3.1]{DiS}
\label{def:DS:3.1}
An infinite cardinal $\tau$ is called \emph{admissible\/} if there 
	exists a pseudocompact group of cardinality $\tau.$
\end{definition}

\begin{definition}\cite{{CEG}, DiS}
Let $\sigma$ be an infinite cardinal.
	The condition $Ps(\tau,\sigma)$ states that the
	set $\{0,1\}^\sigma$ contains a subset of cardinality
	$\tau$ whose projection on every countable subproduct
	$\{0,1\}^A$ is surjective. Furthermore, we define
$$
\delta(\sigma)=\min\{\tau: Ps(\tau,\sigma)
\text{ holds}\}.
$$	
\end{definition}

\begin{theorem}
\label{3,2.7}
Let $\tau$ be an infinite cardinal.
\begin{itemize}
\item[(i)] \cite[Lemma 2.7]{DiS}
	If $\sigma\ge\omega$, then
	$Ps(\tau,\sigma )$ is equivalent to $\delta(\sigma)
	\leq\tau\leq 2^\sigma.$
\item[(ii)] \cite[Theorem 3.3]{DiS}
$\tau$ is admissible if and only if 
	the condition  $Ps(\tau,\sigma)$ holds for some
	cardinal
	$\sigma\ge\omega.$ 
\end{itemize}
\end{theorem}

\begin{fact}
\label{nG} 
\label{folklore} 
(i) \cite[Theorem 3.5]{DiS} Let $G$ be a pseudocompact abelian group. Then for every $n\in\N$, the subgroup $nG$ of $G$ is pseudocompact (in the subgroup topology).

(ii) \cite[Corollary 3.6]{DiS} 
If an abelian group $G$ admits a pseudocompact group topology, then so does its subgroup $nG$ for every $n\in\N$. In particular, the cardinal $|nG|$ is either finite or admissible for every $n\in\N$.
\end{fact}
\begin{proof}
(i) The map $g\mapsto ng$ is a continuous surjection from $G$ onto $nG$.
It remains to note that pseudocompactness is preserved by continuous surjections.

(ii) follows from (i) and Definition \ref{def:DS:3.1}.
\end{proof}

\begin{theorem}\cite[7.4]{CoRo}
\label{van:Douwen:theorem}
Every pseudocompact torsion abelian group is bounded.
\end{theorem}

\begin{definition}
\label{UK:definition}
Recall that for each bounded torsion abelian group $B$, there exist a finite set $P_B\subseteq \Prm$, and for every $p\in P_B$, a positive integer $r_p^B$ and  (possibly zero) cardinals 
$\{\alpha_{k,p}^B: k=1,\dots, r_p^B\}$ called {\em Ulm-Kaplanski ivariants\/} of $B_p$ (and of $B$)
such that 
\begin{equation}
\label{eq:18:B}
B=\bigoplus_{p\in P_B} B_p
\ \ \text{ and }\ \ 
B_p=\bigoplus_{i=1}^{r^B_p} \Z(p^i)^{(\alpha_{i,p}^B)}.
\end{equation}
The decomposition \eqref{eq:18:B} of $B$ is unique up to isomorphism, and all Ulm-Kapalanski invariants $\alpha_{i,p}^B$ of $B$ are uniquely determined by $B$.
\end{definition}

\begin{theorem}\label{3,6.2}\cite[Theorem 6.2]{DiS}
	For a bounded torsion abelian group $G$ the following 
	conditions are equivalent:
	\begin{enumerate}
		\item $G$ has a pseudocompact group topology,
		\item for all $n\in\mathbb{N}$, the cardinal $|nG|=|G/G[n]|$ 
		is either finite or admissible,
		\item 
		for each $p\in \mathbb{P}$, the 
		$p$-subgroup 
		$G_p=\bigcup_{k\in\N} G[p^k]$ 
		of $G$ has a pseudocompact group topology, 	
		\item for $k\in\mathbb{N}$ and $p\in\mathbb{P}$, the cardinal $\beta_{k,p}^G=|G_p/G_p[p^k]|=|p^kG_p|$
		is either finite or admissible. 
		When $\beta_{k,p}^G$ is infinite,
		$\beta_{k,p}^G=\max\{\alpha_{j,p}^G: j=k+1,\dots,   r_{p}^G\}$ holds.
	\end{enumerate}
\end{theorem}

An alternative description of bounded torsion groups admitting a pseudocompact group topology was obtained by Comfort and Remus \cite{CR2}.

From Theorems \ref{van:Douwen:theorem} and \ref{3,6.2}, we get the following corollary:
\begin{corollary} 
\label{necessary:and:sufficient:condition:for:torsion:groups}
An abelian torsion group $G$ admits a pseudocompact group topology if and only if $G$ is bounded and the cardinal $|nG|$ is either finite or admissible for every $n\in\N$. 
\end{corollary}

Abelian groups supporting a pseudocompact {\em connected\/} (Hausdorff) group topology 
are completely described in the following theorem:

\begin{theorem}\label{3,7.1}\cite[Theorem 7.1]{DiS}
	For an abelian group $G$ and a cardinal $\sigma\geq \omega_1$, the following conditions are equivalent:
	\begin{enumerate}
		\item $G$ admits a pseudocompact connected group 
		topology of weight $\sigma,$ 
		\item $G$ admits a pseudocompact connected, locally connected group 
		topology of weight $\sigma,$ 
		\item $Ps(r_0(G),\sigma )$ and $|G|\le 2^\sigma$,
		\item $\delta(\sigma)\leq r_0(G)\leq|G|\leq 2^\sigma$. 
	\end{enumerate}
\end{theorem}

By Theorem \ref{3,6.2} and Corollary \ref{necessary:and:sufficient:condition:for:torsion:groups},
the algebraic structure of abelian torsion groups 
having 
pseudocompact group topologies 
is completely clear.
An important necessary condition for the existence of pseudocompact group topologies on 
non-torsion abelian groups
was found by Dikranjan and Giordano Bruno.

\begin{theorem}\cite[Theorem 1.18]{DiG}\label{2,1.18}
	If $G$ is a pseudocompact non-torsion abelian group,
	then $Ps(r_0(G),w_d(G))$ holds, where $r_0(G)$ is the free rank of $G.$ 
	(We refer the reader to Definition \ref{2,1.7} below for notation $w_d(G)$.)
\end{theorem}

A different proof of Theorem \ref{2,1.18} was given 5 years earlier by Dikranjan and Shakhmatov \cite{DiS-preprint}, but it remains unpublished to this day.

For torsion and divisible groups,
``algebra forces topology''.

\begin{theorem}
\label{psc:bounded:zero-dimensional}
\label{Wilcox}
\begin{itemize}
\item[(i)] Every pseudocompact torsion group is zero-dimen\-sional \cite[Theorem 1.8]{Dik}.
\item[(ii)] A pseudocompact divisible group is always connected \cite{Wilcox}.
\end{itemize}
\end{theorem}

However, this rigidity does not always happen. 
\begin{example}
\label{free:group:example}
The free abelian group of the cardinality continuum admits both a zero-dimensional pseudocompact group topology and a connected pseudocompact group topology \cite[Theorem 5.6]{DiS}.
\end{example}

Combining Theorems \ref{3,7.1} and \ref{Wilcox}~(ii), we obtain a characterization of divisible groups having a pseudocompact group topology:

\begin{theorem}\cite[Corollary 7.9]{DiS}
\label{divisible:pseudocompact;groups}
A non-trivial divisible abelian group admits a pseudocompact group topology 
if and only if $\delta(\sigma)\leq r_0(G)\leq|G|\leq 2^\sigma$ for some cardinal $\sigma\ge\omega_1$.
\end{theorem}

The paper is organized as follows. Section 
\ref{w-divisible:section} collects necessary background on \cd\ and $w$-divisible abelian groups.
This section can be lasily skipped at first reading. 
Our principal results are exposed
in Sections \ref{sec:psc:topologies:on:cd:groups},
\ref{sec:complete:solution:via:GM:decompositions}
and
\ref{sec:intrinsic:characterization}.
Corollary \ref{structure:of:card-divisible:groups}
extends Theorem \ref{divisible:pseudocompact;groups} to a much wider class of \cd\ abelian groups (defined in Definition \ref{def:w-divisible}).
Section \ref{sec:GM:decompositions}
introduces direct sum decompositions $G=B\oplus H$ of an abelian group $G$ into a bounded torsion group $B$ and a \cd\ group $H$ (we refer to such a decomposition as a {\em Galindo-Macario decomposition\/} of $G$).

The main result of Section
\ref{sec:complete:solution:via:GM:decompositions}
is Corollary
\ref{equivalent:conditions:for:psc:topologizations}
giving a complete desctiption of abelian groups
$G$ which can be equipped with a pseudocompact group topology in terms of 
Galindo-Macario decompositions of $G$,
thereby resolving 
Problem \ref{main:problem}. As witnessed by 
Theorem \ref{thm:6.2},
a crucial role here is played by what we call
{\em minimal\/} Galindo-Macario decompositions, i.e.,
Galindo-Macario decompositions minimizing the order of the bounded part $B$ of the decomposition (Definition \ref{def:b_G});
we shall denote that minimal order by $b_G$.

In Section \ref{sec:intrinsic:characterization} we give
an ``intrinsic'' characterization of abelian groups $G$ admitting a pseudocompact group topology in terms of a {\em minimal\/}
positive integer $m_G$ such that $m_G G$ is \cd. It turns out that an abelian group $G$ can be equipped with a pseudocompact group topology iff both its bounded subgroup $G[m_G]$ and its \cd\ subgroup $m_G G\cong G/G[m_G]$ admit a pseudocompact group topology (Theorem \ref{no:GM}). Item (ii) of Corollary \ref{internal:description} offers equivalent algebraic and set-theoretic conditions
for this to happen. 

In Section \ref{minimal:GM:decompositions}
we 
study minimal Galindo-Macario decompositions of 
abelian groups.
Theorem \ref{characterizations:of:minimality}
offers a 
couple of
equivalent conditions to the minimality of a Galindo-Macario decomposition
of an 
abelian 
group.
These characterizations are essentially used in 
the 
proofs of 
Theorems 
\ref{thm:6.2}
and \ref{no:GM}
in Section
\ref{sec:second:group:of:proofs}.
Corollary \ref{algebraic:corollary} 
establishes
the equality $b_G=m_G$ for each abelian group $G$.

\section{$\varphi$-divisible, \cd\ and $w$-divisible groups}
\label{w-divisible:section}
\label{Sec.5}

Let $\varphi$ be a cardinal function defined on the class of all abelian topological groups which is {\em monotone\/} in a sense that $\varphi(H)\le \varphi(G)$ whenever $H$ is a subgroup of $G$. 
\begin{definition}
\label{def:phi_d}
Following \cite[Remark 7.9]{DS-Advances}, for every topological group $G$, we can define the cardinal 
\begin{equation}
\label{eq:1:fj}
\varphi_d(G)=\min\{\varphi(nG):n\in\N\}.
\end{equation}
\end{definition}
One can consider the new cradinal function $\varphi_d$ to be a ``divisible modification'' of the original cardinal function $\varphi$.

It is clear from \eqref{eq:1:fj} that $\varphi_d(G)\le \varphi(G)$ for every topological group $G$.
\begin{definition}
\label{def:phi-divisible}
Following \cite[Remark 7.9]{DS-Advances}, we say that a topological group $G$ is {\em $\varphi$-divisible\/} provided that $\varphi_d(G)=\varphi(G)$, or equivalently, if 
$\varphi(nG)=\varphi(G)$ for every $n\in\N$. 
\end{definition}
Since $nG=G$ for a divisible group $G$, 
all divisible groups are $\varphi$-divisible, thereby justifying the appearence of the word ``divisible'' in the notion of $\varphi$-divisibility.

\begin{lemma}
\label{general:lemma:on:multiples:new}
Let $\varphi$ be a monotone cardinal function defined on the class of topological abelian groups and let $G$ be an abelian topological group.
Then the following holds:
\begin{itemize}
\item[(i)] if $H$ is a subgroup of $G$, then 
$\varphi_d(H)\le \varphi_d(G)$; that is,
$\varphi_d$ is a monotone cardinal function on the class of topological abelian groups.
\item[(ii)]
$\varphi_d(kG)=\varphi_d(G)$ for every $k\in\N$.
\item[(iii)]
if $n\in\N$ and $\varphi_d(G)=\varphi(nG)$, then the subgroup $nG$ of $G$ is $\varphi$-divisible  (in its subgroup topology).
\item[(iv)] if $G$ is $\varphi$-divisible, then so is its subgroup $kG$ for every $k\in\N$.
\end{itemize}
\end{lemma}
\begin{proof}
(i) Let $n\in\N$ be arbitrary. Then $nH$ is a subgroup of $nG$, so $\varphi(nH)\le \varphi(nG)$ by monotonicity of $\varphi$. Therefore,
$$
\varphi_d(H)=\min\{\varphi(nH):n\in\N\}\le \min\{\varphi(nG):n\in\N\}=\varphi_d(G)
$$
by \eqref{eq:1:fj}.

(ii)
Fix $k\in\N$. 
Then 
$kG$
is a subgroup of 
$G$, so 
$\varphi_d(kG)\le \varphi_d(G)$ by (i).
In the opposite direction, the inequality
$$
\varphi_d(G)=\min\{\varphi(nG):n\in\N\}\le \min\{\varphi(nkG):n\in\N\}=\varphi_d(kG)
$$
is clear from the definition of the function minimum and \eqref{eq:1:fj}.

(iii) Suppose that $n\in\N$ and $\varphi_d(G)=\varphi(nG)$. Then $\varphi_d(nG)=\varphi_d(G)=\varphi(nG)$ by (i) and our assumption. This shows that 
$nG$ is $\varphi$-divisible by Definition \ref{def:phi-divisible}.

(iv) Suppose that $G$ is $\varphi$-divisible and $k\in\N$. Then $\varphi_d(kG)\le \varphi(kG)$ by \eqref{eq:1:fj}, and $\varphi_d(G)=\varphi_d(kG)$ by (i). Since $G$ is $\varphi$-divisible, $\varphi(G)=\varphi_d(G)$. Finally, since $kG$ is a subgroup of $G$ and $\varphi$ is monotone, $\varphi(kG)\le \varphi(G)$. Combining all these inequalities, we deduce that $\varphi_d(kG)=\varphi(kG)$. Therefore, $kG$ is $\varphi$-divisible.
\end{proof}

From 
\eqref{eq:1:fj}
and Lemma \ref{general:lemma:on:multiples:new}(iii) we get the following
\begin{corollary}
\label{cor:phi-divisible-multiples}
Let $\varphi$ be a monotone cardinal finction defined on the class of topological abelian groups. Then for every  abelian topological group $G$ there exists $n\in\N$ such that the subgroup $nG$ of $G$ is $\varphi$-divisible  (in its subgroup topology).
\end{corollary}

For a topological abelian group $G$, let 
$w(G)$  denote the weight of $G$.
Applying Definitions \ref{def:phi_d} and \ref{def:phi-divisible} to the monotone cardinal function $\varphi$ defined by $\varphi(G)=w(G)$, 
we get the important particular case introduced in \cite{DiG}.
\begin{definition}
\cite{DiG}
\label{2,1.7}
Let $G$ be a topological abelian group.
The cardinal
$w_d(G)=\min\{w(nG):n\in\mathbb{N}\}$
is called the {\em divisible weight\/} of$G$.
If $w_d(G)=w(G)$, then $G$ is called
{\em $w$-divisible\/}.
\end{definition}
\begin{remark}
In \cite[Definition 1.7]{DiG},
the authors 
introduced
the divisible weight of $G$ by 
$w_d(G)=\min\{w(m!G):m\in\N\}$. Let us check that 
$\min\{w(nG):n\in\mathbb{N}\}=\min\{w(m!G):m\in\N\}$, so these two definitions are equivalent.
The inequality $\min\{w(nG):n\in\mathbb{N}\}\le\min\{w(m!G):m\in\N\}$ is trivial.
For every $n\in \N$, note that $n!G$ is a subgroup of $nG$, so $w(n!G)\le w(nG)$.
Therefore, $\min\{w(n!G):n\in\N\}\le \min\{w(nG):n\in\N\}$.
\end{remark}

Every abelian group 
can be viewed as a topological group when equipped with its discrete 
topology. We shall use this identification, thereby making 
Definitions \ref{def:phi_d} and \ref{def:phi-divisible}
applicable also
to abstract abelian groups.

\begin{definition}
\label{def:w-divisible}
Let $G$ be an abelian (topological) group. 
For the cardinal
function $\varphi$ defined by $\varphi(G)=|G|$, we denote the value of its ``divisible modification'' $\varphi_d$ on 
$G$ by $|G|_d$; that is, we let
\begin{equation}
|G|_d=\min\{|nG|: n\in\N\}.
\end{equation}
We say that $G$ is {\em \cd\/} if $|G|_d=|G|$, or equivalently, 
	$|nG|=|G|$ for every $n\in\mathbb{N}$.
\end{definition}

Applying Corollary \ref{cor:phi-divisible-multiples} to 
this function $\varphi$ from Definition \ref{def:w-divisible},
we obtain 
the following 
\begin{corollary}
\label{cd:core}
For every  abelian group $G$, there exists $n\in\N$ such that the subgroup $nG$ of $G$ is \cd.
\end{corollary}

Let $G$ be an 
abelian group equipped with its discrete topology. Then
$w(G)=|G|$, so $w_d(G)=|G|_d$. Therefore, our notion of divisible cardinality 
$|G|_d$  
for an abstract abelian group $G$ is a particular (discrete) case of the notion of divisible weight $w_d(G)$ 
introduced in \cite{DiG}. Unfortunately, the same notation $w_d(G)$ is used in the literature for denoting the divisible weight of both topological groups and abstract (discrete) groups. When $G$ is a topological group, this creates an ambiguity, as 
$|G|_d$
and $w_d(G)$ need not coincide for $G$.
So far this never created a problem beause in each particular situation only one version (either topological $w_d(G)$ or algebraic $|G|_d$) was used at a time.
However, in the present paper (to be precise, in the proof of Theorem \ref{clever:lemma}) we need to use {\em both\/} of these notions at the same time. 
This justifies the introduction of a new notation $|G|_d$ in Definition \ref{def:w-divisible}.

\begin{remark}
\label{no:finite:divisble:weights}
\label{countable:cd:trivial}
\emph{(i) For every abelian group $G$, either $|G|_d=1$ (and $G$ is a bounded torsion group) or $|G|_d\ge\omega$.}
Indeed, if $|G|_d>1$, then $|G|_d$ cannot be finite, as every finite group is bounded torsion.
In particular, {\em $|G|_d\ge \omega$ for every unbounded abelian group $G$.}

(ii) {\em A bounded \cd\ abelian group $H$ is trivial.\/}
Indeed, $|H|_d=1$ by (i), as $H$ is bounded.
Since $H$ is \cd, $|H|=|H|_d$ by Definition \ref{def:w-divisible}.
\end{remark}

\begin{lemma}
\label{adding:small:piece:to:w^divisible:group}
If $H$ is a 
\cd\ 
abelian group and $G$ is an abelian group satisfying 
$|G|\le |H|$, then the direct sum $G\oplus H$ is also 
\cd. 
\end{lemma}
\begin{proof}
Our lemma obviously holds when $H$ is trivial, so we shall assume now that $H$ is non-trivial.
Since $H$ is non-trivial and 
\cd, 
it must be infinite by Remark \ref{countable:cd:trivial}~(ii).

Let $n\in \N$. 
Since $H$ is \cd.
$|nH|=|H|\ge \omega$ by Definition \ref{def:w-divisible}.
Note that $|nG|\le |G|\le |H|$ and $n(G\oplus H)=nG\oplus nH$, so
$
|n(G\oplus H)|=\max\{|nG|, |nH|\}=\max\{|nG|, |H|\}=|H|=|G\oplus H|.
$
Therefore, $G\oplus H$ is 
\cd\ 
by Definition \ref{def:w-divisible}.
\end{proof}

When $r(G)$ denotes the rank of an abelian group $G$, its {\em divisible rank\/} $r_d(G)$
was introduced in
\cite{Szele} (in the case of discrete $p$-groups)  under the name of {\em final rank\/} and later studied in some detail in 
\cite[Section 7]{DS-Advances}. 

\section{Pseudocompact group topologies on \cd\ abelian groups}
\label{sec:psc:topologies:on:cd:groups}

We refer to Definitions \ref{2,1.7} and 
\ref{def:w-divisible} for the notions appearing in the next theorem.
\begin{theorem}
\label{clever:lemma}
Assume that $m\in\N$, $H$ is 
a non-trivial
\cd\ abelian 
group such that 
its subgroup $mH$ admits a 
pseudocompact group topology of divisible weight $\mu$.
Then $H$ 
can be equipped with a connected, locally connected pseudocompact group topology of weight 
\begin{equation}
\label{eq:15:fd}
\sigma=\max\{\mu,\omega_1\}.
\end{equation}
\end{theorem}
\begin{proof}
Since a bounded \cd\ abelian group is trivial by Remark \ref{countable:cd:trivial}~(ii),
$H$ must be unbounded.
Thus, $mH$ is unbounded as well.
Fix a 
pseudocompact group topology of divisible weight $\mu$ on $mH$.
Then $mH$ is non-torsion by 
Theorem \ref{van:Douwen:theorem}.
Applying Theorem \ref{2,1.18} to the 
non-torsion pseudocompact 
group $mH$,
we conclude that
$$
Ps(r_0(mH),w_d(mH))=Ps(r_0(mH),\mu)
$$
holds.
By Theorem \ref{3,2.7}~(i), this implies 
\begin{equation}
\label{17:sd}
\delta(\mu)\le r_0(mH).
\end{equation}
Since $mH$ is a subgroup of $H$, it follows that 
\begin{equation}
\label{19:sd}
r_0(mH)\le r_0(H)\le |H|.
\end{equation}

By Definition \ref{2,1.7}, 
there exists $n\in\N$ such that 
$w(nmH)=w_d(mH)=\mu$.
Since $w(nmH)=\mu$,
the cardinality of the subgroup $nmH$ of $H$ 
satisfies 
\begin{equation}
\label{18:sd}
|nmH|\le 2^\mu.
\end{equation}
Since $H$ is
\cd\ 
by the assumption of our lemma,
we have 
\begin{equation}
\label{20:sd}
|H|=|nmH|.
\end{equation}
Combining \eqref{17:sd}, \eqref{19:sd},  \eqref{20:sd} and \eqref{18:sd}, we get 
\begin{equation}
\label{key:inequalities:on:H:mu}
\delta(\mu)\leq r_0(H)\leq |H|
\leq 2^{\mu}.
\end{equation}
\begin{claim}
\label{key:inequalities:on:H:sigma}
$\delta(\sigma)\leq r_0(H)\leq |H|
\leq 2^{\sigma}$.
\end{claim}
\begin{proof}
If $\mu\ge \omega_1$, then $\sigma=\mu$ by \eqref{eq:15:fd}, so 
our claim 
coincides with \eqref{key:inequalities:on:H:mu}.
Suppose 
that $\mu=\omega$.
In this case $\sigma=\omega_1$ by \eqref{eq:15:fd},
so
$\delta(\mu)=\delta(\omega)=\cont=\delta(\omega_1)=\delta(\sigma)$
by \cite[Lemmas 2.7 and 2.8]{DiS}.
Since
$2^\mu=2^\omega\le 2^{\omega_1}=2^\sigma$,
our claim
follows from \eqref{key:inequalities:on:H:mu}.
\end{proof}

Note that $\sigma\ge \omega_1$ by \eqref{eq:15:fd}.
 By Claim \ref{key:inequalities:on:H:sigma} and the implication (3)$\to$(2) of Theorem \ref{3,7.1} applied to $H$, we conclude that $H$ 
can be equipped with 
 a pseudocompact connected, locally connected group topology of weight $\sigma$.
\end{proof}

The following corollary 
extends Theorem \ref{divisible:pseudocompact;groups}
from divisible groups to the much wider class of
\cd\ groups.  
\begin{corollary}
\label{structure:of:card-divisible:groups}
let $G$ be a 
non-trivial 
\cd\ abelian group. Then the following conditions are equivalent:
\begin{itemize}
\item[(i)] $G$ admits a pseudocompact group topology;
\item[(ii)] $G$ admits a pseudocompact connected, locally connected group topology of uncountable weight;
\item[(iii)] $\delta(\sigma)\leq r_0(G)\leq|G|\leq 2^\sigma$ for some cardinal $\sigma\ge \omega_1$.
\end{itemize}
\end{corollary}
\begin{proof}
The
implication (i)$\to$(ii) follows from Theorem \ref{clever:lemma} (in which one takes $m=1$).
The implications (ii)$\to$(iii)$\to$(i) follow from 
Theorem \ref{3,7.1}.
\end{proof}

\begin{remark}
Since the free abelian group 
of cardinality continuum is clearly \cd, Example \ref{free:group:example} shows that Theorem \ref{Wilcox}~(ii) cannot be extended from divisible groups to \cd\ groups. However,
Corollary \ref{structure:of:card-divisible:groups} shows that some weaker property is retained by all 
\cd\  abelian groups; namely, if a non-trivial \cd\ abelian groups has a pseudocompact group topology, it also has (possibly) another topology which is connected (and locally connected).
\end{remark}

\section{Galindo-Macario decompositions and their minimality}
\label{sec:GM:decompositions}

\begin{definition}
\label{def:GM:decomposition}
Let $G$ be an abelian group.
A direct sum decomposition 
\begin{equation}
\label{eq:9:si}
G=B\oplus H
\end{equation}
into a bounded torsion group $B$ and a \cd\ group $H$ will be referred to as a {\em Galindo-Macario decomposition\/} of $G$. 
\end{definition}

The reason for the terminology in Definition
\ref{def:GM:decomposition} is rooted in the following result
discovered by 
Galindo and Macario (who clearly did not use our terminology):
\begin{fact}
\label{Galindo-Macario-theorem}
\cite[Lemma 4.1]{GM}
Every abelian group $G$ 
has a Galindo-Macario decomposition.
\end{fact}

An alternative proof of this Fact (and a ``dual version'' of it for compact groups) can be found in \cite{DiG2}.
Fact \ref{Galindo-Macario-theorem} has already found its applications in 
\cite{DS-Advances} and \cite{DS-Minap}. 

A
Galindo-Macario decomposition from 
Definition
\ref{def:GM:decomposition}
is far from unique.

\begin{lemma}
\label{powers:of:GM:decompositions} 
\label{rem:cardinality:of:H} 
Let
$G=B\oplus H$ be a Galindo-Macario decomposition of an abelian group $G$. Then the following holds:
\begin{itemize}
\item[(i)]
If $n\in \N$ and $nB=\{0\}$,
then $nG=nH$.
\item[(ii)] $|H|=|G|_d$.
\end{itemize}
\end{lemma}
\begin{proof}
(i)
Follows 
from \eqref{eq:9:si} and $nB=\{0\}$.

\smallskip
(ii)
Since
$B$ is bounded, $nB=\{0\}$ for some $n\in\N$. 
Then $nG=nH$ by (i),
so $|G|_d=|nG|_d=|nH|_d=|H|_d$ by Lemma \ref{general:lemma:on:multiples:new}(ii) applied twice to the monotone cardinal function 
$\varphi(G)=|G|$. Since $H$ is \cd, $|H|_d=|H|$ by 
Definition \ref{def:w-divisible}.
This proves the equality $|G|_d=|H|$.
\end{proof}

\begin{definition}
\label{def:b_G}
For an abelian group $G$,
we let 
\begin{equation}
\label{eq:33:fo}
b_G=\min\{
n_B:
G=B\oplus H \text{ is Galindo-Macario decomposition of } G\},
\end{equation}
where $n_B$ denotes the order of $B$.
A Galindo-Macario decomposition $G=B\oplus H$
of $G$ satisfying $n_B=b_G$ will be called {\em minimal\/}.
\end{definition}

\begin{remark}
\label{rem:minimal:decomposition}
Let $G$ be an abelian group.

(i) It follows from Fact \ref{Galindo-Macario-theorem} that $G$ has at least one Galindo-Macario decomposition. This means that {\em the integer 
$b_G$ from 
\eqref{eq:33:fo} is well-defined\/}.

(ii) 
{\em 
$G$ has a minimal Galindo-Macario decomposition\/}. Indeed, it follows from 
\eqref{eq:33:fo} that 
there exists a Galindo-Macario decomposition $G=B\oplus H$
of $G$ satisfying $n_B=b_G$. Now the decomposition 
$G=B\oplus H$ is minimal by Definition \ref{def:b_G}.
\end{remark}

\section{Solving Problem \ref{main:problem}
via Galindo-Macario decompositions}
\label{sec:complete:solution:via:GM:decompositions}

In
this 
section
we
study topological properties of the summands $B$ and $H$ of Galindo-Macario decompositions $G=B\oplus H$ of a pseudocompact abelian group $G$.

Our first decomposition result shows that, somewhat surprisingly, for {\em each\/} Galindo-Macario decomposition
$G=B\oplus H$ of 
a 
pseudocompact abelian group $G$, its \cd\ part $H$ {\em automatically\/} admits a pseudocompact group topology.

\begin{theorem}\label{weight:divisible:part}
Let $G=B\oplus H$ be a Galindo-Macario decomposition of a pseudocompact abelian group $G$.
Then $H$ admits a pseudocompact connected, locally connected group topology.
\end{theorem}
\begin{proof}
If $H$ is trivial, then the only (group) topology on $H$ is compact (so pseudocompact by \eqref{compact:implies:pseudocompact}), connected and locally connected.
Suppose now that $H$ is non-trivial.
Since $B$ is bounded,
there exists 
$m\in\mathbb{N}$ such that $mB=\{0\}$, 
so $mG=mH$ by Lemma \ref{powers:of:GM:decompositions}~(i).
Since $G$ is pseudocompact, 
the subgroup $mG=mH$ of $G$ is pseudocompact (in its subgroup topology) by
Fact \ref{nG}~(i). 
Applying Theorem \ref{clever:lemma}, we conclude that $H$ 
 admits a connected, locally connected pseudocompact group topology.
\end{proof}
Our second decomposition result shows that, for a pseudocompact abelian group $G$, 
one can {\em choose\/} a Galindo-Macario decomposition 
in such a way that 
its bounded part
$B$ 
admits
a pseudocompact group topology.

\begin{theorem}
\label{thm:6.2}
Let $G=B\oplus H$ be a minimal Galindo-Macario decomposition of a pseudocompact abelian group $G$.
Then $B$ admits a (zero-dimensional)
pseudocompact group topology.  
\end{theorem}

The 
proof
of 
Theorem 
\ref{thm:6.2}
is
postponed until Section 
\ref{sec:second:group:of:proofs}.

\begin{definition}
A subgroup $H$ of a topological group $G$ {\em topologically splits in $G$\/} if 
there exists a topological group $B$ 
such that $G$ is topologically ispomorphic to the direct product $B\times H$ equipped with the Tychonoff product topology.
\end{definition}

\begin{proposition}
\label{direct:product:proposition}
Suppose that an abelian group $G$ is decomposed as in \eqref{eq:9:si},
where 
$B$ 
admits a pseudocompact group topology $\mathscr{T}_B$ and 
$H$ 
admits a 
pseudocompact group topology $\mathscr{T}_H$.
Note that $G$ can be viewed as a direct product $B\times H$ and be equipped with the product topology 
$\mathscr{T}_G$ of the direct product $(B,\mathscr{T}_B)\times (H,\mathscr{T}_H)$.
Then
$\mathscr{T}_G$ is a pseudocompact group topology on $G$ such that both $(B,\mathscr{T}_B)$ and $(H,\mathscr{T}_H)$ topologically split in $(G,\mathscr{T}_G)$.
Furthermore, if $B$ is a bounded (torsion) group
and the topology
$\mathscr{T}_H$ is connected,
then 
$(H,\mathscr{T}_H)$ 
coincides with the connected component of $0$ in $(G,\mathscr{T}_G)$. 
\end{proposition}
\begin{proof}
Since both $\mathscr{T}_B$ and $\mathscr{T}_H$ are pseudocompact, 
$\mathscr{T}_G$ is pseudocompact by 
Theorem \ref{Comfort-Ross-theorem}.

Suppose now that
$B$ is a bounded (torsion) group
and the topology
$\mathscr{T}_H$ is connected.
Since $B$ is a bounded torsion group, the pseudocompact group topology $\mathscr{T}_B$ is zero-dimensional by Theorem \ref{psc:bounded:zero-dimensional}~(i). 
Therefore, $(H,\mathscr{T}_H)$ is the maximal connected subgroup of $(G,\mathscr{T}_G)$ containing $0$.
\end{proof}

\begin{corollary}
\label{equivalent:conditions:for:psc:topologizations}
For an abelian group $G$, the following 
conditions 
are equivalent:
\begin{itemize}
\item[(i)] $G$ admits a pseudocompact group topology;
\item[(ii)] $G$ can be decomposed into a direct sum \eqref{eq:9:si}
of 
	a  bounded
	torsion group $B$ and a \cd\
	group $H$ which both admit a pseudocompact group topology;
\item[(iii)] $G$ can be decomposed into a direct sum \eqref{eq:9:si},
	where $B$ is 
	a  bounded
	torsion group admitting a (zero-dimensional) pseudocompact 
group topology and $H$ is
	a \cd\ 
	group admitting a connected, locally connected pseudocompact group topology;
\item[(iv)] $G$ can be equipped with a pseudocompact group topology $\mathscr{T}$ such that the connected component 
of $(G,\mathscr{T})$  is locally connected, \cd\ 
and topologically splits in $(G,\mathscr{T})$.
\end{itemize}
\end{corollary}
\begin{proof}
(i)$\to$(iii) Let $G=B\oplus H$ be a minimal Galindo-Macario decomposition of $G$ which exists by Remark \ref{rem:minimal:decomposition}~(ii).
Since $G$ admits a pseudocompact group topology by (i), $B$ admits a (zero-dimensional) pseudocompact 
group topology by Theorem \ref{thm:6.2}, and $H$ admits a connected, locally connected pseudocompact group topology by Theorem
\ref{weight:divisible:part}.

The implications (ii)$\to$(i) 
and
(iii)$\to$(iv) follow from Proposition \ref{direct:product:proposition},
and the implications
(iii)$\to$(ii)
and
(iv)$\to$(i) are trivial.
\end{proof}

\begin{remark}
\label{long:remark}
(i)
{\em If the connected component $c(G)$ of a pseudocompact group $G$ topologically splits in $G$,
then $c(G)$ must be pseudocompact as well.\/} Indeed, in this case $c(G)$ is also a quotient of $G$, and pseudocompactness is preserved by quotients.

(ii)
The  pseudocompact group topologies provided on an abelian group $G$
by item (iv) of
Corollary 
\ref{equivalent:conditions:for:psc:topologizations}
have an additional property that their connected component $c(G)$ 
topologically splits 
in
$G$, so $c(G)$ is itself pseudocompact by (i).
This does not hold in general.

(iii)
It was proved in \cite[Corollary 3.6]{Scorza} (announced also in \cite[Corollary 1.4.6]{D_CC}) that
for every precompact connected abelian group  $L$ and for every uncountable cardinal $\sigma \geq w(L)$,
there exists a pseudocompact abelian group $G$  of weight $\sigma$  such that $c(G)=q(G) \cong L$.
(The particular (weight-free) case of this result for connected groups $L$ such that the completion $\widetilde L$ is torsion-free  was proved in \cite[Theorem 7.6]{CvM2}, reinforcing \cite[Lemma 2.1]{CvM1}, claiming that every precompact group can be presented as a
closed subgroup of a pseudocompact group.)

(iv) It follows from (iii) that {\em 
the connected component $c(G)$ of a pseudocompact abelian group $G$ need not be pseudocompact, and thus need not topologically split in $G$ by (i).\/}
\end{remark}

\begin{remark}
When combined with  Theorems \ref{3,6.2} and  \ref{3,7.1} (or Corollary \ref{structure:of:card-divisible:groups}),
our Corollary \ref{equivalent:conditions:for:psc:topologizations}
offers a complete description of the algebraic structure of pseudocompact abelian groups, thereby solving Problem \ref{main:problem}.
A different
description of the algebraic structire of pseudocompact abelian groups was obtained by Dikranjan and Shakhmatov in their preprint
\cite{DiS-preprint} (cited as [20] in \cite{DiG} and as [14] in \cite{GM}) but it remained unpublished to this day. 
\end{remark}

\section{An integer 
determining
the fate of Problem \ref{main:problem}}
\label{sec:intrinsic:characterization}

In this section we give an ``intrinsic'' characterization of abelian groups admitting a pasedocompact group topology. A principal role in this characterization is played by the integer number introduced  next.
\begin{definition}
\label{definition:of:m_G}
For every abelian group $G$, there exists $n\in\N$ such that $nG$ is \cd\ (see Corollary \ref{cd:core} below), so the positive integer 
\begin{equation}
\label{eq:1:fo}
m_G=\min\{n\in\N: nG \text{ is \cd}\}
\end{equation}
is well-defined.
\end{definition}

The implication (ii)$\to$(i) of our next theorem provides a positive answer to \cite[Problem 9.10]{DiS}
for $n=m_G$.

\begin{theorem}
\label{no:GM}
Let $G$ be an abelian group, and 
let $m=m_G$ for breivity.
Then the following conditions are equivalent:
\begin{itemize}
\item[(i)]
$G$ admits a pseudocompact group topology;
\item[(ii)]
both the bounded group $G[m]$ and the \cd\ group
$G/G[m]\cong mG$
admit a pseudocompact group topology.
\item[(iii)] $G[m]$ has a zero-dimensional pseudocompact group topology and $md\cong G/G[m]$ has a connected, locally connected pseudocompact group topology.
\end{itemize}
\end{theorem}

The proof of this theorem is postponed until Section 
\ref{sec:second:group:of:proofs}.

\begin{corollary}
\label{internal:description}
Let $G$ be
abelian group, and let $m=m_G$ for breivity.
Then the following conditions are equivalent:
\begin{itemize}
\item[(i)]
$G$ admits a pseudocompact group topology;
\item[(ii)] $G$ satisfies the following two conditions:
\begin{itemize}
\item[(ii-a)] if $mG\neq\{0\}$, then $\delta(\sigma)\leq r_0(G)\le|mG|\leq 2^\sigma$ for some cardinal $\sigma\ge \omega_1$, and 
\item[(ii-b)] 
for each $n\in\N$, the cardinal
$|n G[m]|$ is either finite or  
admissible.
\end{itemize}
\end{itemize}
\end{corollary}
\begin{proof}
By Theorem \ref{no:GM}, $G$ admits a pseudocompact group topology if and only if both of the following conditions hold:
\begin{itemize}
\item[(a)] $mG$ admits a pseudocompact group topology, and
\item[(b)] 
$G[m]$ 
admits a pseudocompact group topology.
\end{itemize}

It follows  from our definition of the integer $m$ that $mG$ is \cd. 
If $mG=\{0\}$, then the only group topology on $mG$ is compact, so also pseudocompact by \eqref{compact:implies:pseudocompact}.
Suppose now that 
$mG$ is non-trivial.
By Corollary \ref{structure:of:card-divisible:groups}, $mG$ admits a pseudocompact group topology precisely when 
$\delta(\sigma)\leq r_0(G)\leq|mG|\leq 2^\sigma$ for some cardinal $\sigma\ge \omega_1$.
We have proved that (a) is equivalent to (ii-a).

Since $G[m]$ is a bounded torsion abelian group, 
(b) is equivalent to (ii-b) by Corollary \ref{necessary:and:sufficient:condition:for:torsion:groups}.
\end{proof}

\section{Minimal Galindo-Macario decompositions}
\label{minimal:GM:decompositions}

\def\md{minimal Galinco-Macario decomposition}

In this section we study minimal Galindo-Macario decompositions introduced in  Definition  \ref{def:b_G}. The main result here is Theorem 
\ref{characterizations:of:minimality} below.
\begin{lemma}
\label{m_G:is:a:lower:bound:for:n_B}
For every Galindo-Macario decomposition $G=B\oplus H$ of an abelian group $G$, the inequality $m_G\le n_B$ holds.
\end{lemma}
\begin{proof}
Since $n_B B=\{0\}$, we have $n_B G=n_B H$ by Lemma \ref{powers:of:GM:decompositions}~(ii).
Since $H$ is \cd, so is $n_B H$ by Lemma \ref{general:lemma:on:multiples:new}~(iv) applied to the monotone cardinal function $\varphi(G)=|G|$.
Therefore, $n_B G$ 
is \cd. This implies the inequality $m_G\le n_B$ by 
Definition \ref{definition:of:m_G}.
\end{proof}

\begin{remark}
\label{rem:minimal:decomposition:iii} 
\label{computing:m_G:b_G:for:a:bounded:G} 
{\em A bounded torsion group 
$G$ has a unique Galindo-Macario decomposition (which is necessarily minimal by Definition  \ref{def:b_G}) satisfying $n_B=m_G$\/}.
Indeed, let $G=B\oplus H$ be a Galindo-Macario decomposition of $G$. Then $H$ is both bounded
(as a subgroup of a bounded group $G$) and \cd\ 
(by Definition \ref{def:GM:decomposition}), so 
$H=\{0\}$ by Remark \ref{countable:cd:trivial}~(ii). This implies $B=G$.
Now
$b_G=n_B$ by \eqref{eq:33:fo}.
Since $G$ is bounded, $|G|_d=1$ by Remark
\ref{no:finite:divisble:weights}~(i),
so 
$m_G=\min\{n: nG=\{0\}\}=\min\{n: nB=\{0\}\}=n_B.$
\end{remark}

\begin{theorem}
\label{characterizations:of:minimality}
For
a Galindo-Macario decomposition \eqref{eq:9:si}
of an 
abelian group $G$,
the following conditions are equivalent:
\begin{itemize}
\item[(i)] \eqref{eq:9:si}
is minimal; that is,
$n_B=b_G$ by Definition \ref{def:b_G}.
\item[(ii)] $n_B=m_G$, where
$m_G$ is as in 
Definition \ref{definition:of:m_G}.
\item[(iii)] If $n\in\N$ and $|nB|\le |H|$, then $n_B$ divides $n$; in particular, $nB=\{0\}$.
\end{itemize}
Furthermore, when $G$ is unbounded, then the following item can be added to the list of equivalent conditions (i), (ii) and (iii):
\begin{itemize}
\item[(iv)] all leading Ulm-Kaplanski invariants of $B$ are strictly bigger than $|G|_d=|H|$.
\end{itemize}
\end{theorem}
\begin{proof}
First we consider the special case when $G$ is bounded. By Remark \ref{rem:minimal:decomposition:iii}, $G$ has a unique Galindo-Macario decomposition $G=B\oplus H$
satisfying items (i), (ii) and $H=\{0\}$.
To check (iii),
suppose that $n\in\N$ and $|nB|\le |H|$.
Then $|nB|\le |H|=1$ by our assumption, which implies 
$nB=\{0\}$. Since $n_B$ is the order of $B$, this means that $n_B$ divides $n$.

{\em In the rest of the proof we assume that $G$ is unbounded\/}, so $|G|_d\ge \omega$ by Remark \ref{no:finite:divisble:weights}~(i)
and
$|H|=|G|_d$ by Lemma \ref{rem:cardinality:of:H}~(ii).
Below we use notations of Definition \ref{UK:definition}.

\smallskip
(i)$\to$(iv) 
Let \eqref{eq:18:B}
 be a decomposition of $B$ using the notation of Defintion
\ref{UK:definition}. 
Assume that (iii) fails. Then there exists $q\in P_B$ such that 
$\alpha_{r_q^B,q}^B\le |G|_d$.
Since $H$ is infinite and $\alpha_{r_q^B,q}^B\le |G|_d=|H|$ by our assumption,
$|\mathbb{Z}(q^{r_q^B})^{(\alpha_{r_q^B,q}^G)}|\le |H|$. Since $H$ is \cd,
\begin{equation}
\label{eq:20:H'}
H'=\mathbb{Z}(q^{r_q^B})^{(\alpha_{r_q^B,q}^G)}\oplus H
\end{equation}
is \cd\
by 
Lemma \ref{adding:small:piece:to:w^divisible:group}.
Note that 
\begin{equation}
\label{eq:19:B'}
B'=\left(\bigoplus_{i=1}^{r^B_q-1} \Z(q^i)^{(\alpha_{i,q}^B)}\right)\oplus\left(\bigoplus _{p\in P_B\setminus\{q\}} B_p\right)
\end{equation}
is bounded, as a subgroup of a bounded group $B$.
Clearly, $B=B'\oplus\mathbb{Z}(p^{r_q^B})^{(\alpha_{r_q^B,q}^G)}$.
It follows from \eqref{eq:18:B},
\eqref{eq:19:B'} 
and \eqref{eq:20:H'} that $B'\oplus H'=B\oplus H=G$
is a Galindo-Macario decomposition of $G$ such that
$n_{B'}=n_B/q<n_B$.
Since $n_B=b_G$ by (i), we get $n_{B'}<b_G$,
in contradiction with \eqref{eq:33:fo}.
Therefore, (iv) must hold.

\smallskip
(iv)$\to$(iii)
Assume that $n\in\N$ and $|nB|\le |H|$.
Let $p\in P_B$ be arbitrary.
Since $\Z(p^{r_p^B})^{(\alpha_{{r_p^B},p})}$
is a direct summand of $B$ by 
Definition \ref{UK:definition},
$n \Z(p^{r_p^B})^{(\alpha_{{r_p^B},p})}$
is a direct summand of $nB$. 

Suppose that
$n \Z(p^{r_p^B})^{(\alpha_{{r_p^B},p})}\neq \{0\}$.
Then 
$
\alpha_{r_p^B,p}\le |n\Z(p^{r_p^B})^{(\alpha_{r_p^B},p)}|\le |nB|\le |H|
=|G|_d
$
by our assumption on $n$ and $B$.
On the other hand, $\alpha_{r_p^B,p}>|G|_d$ by 
(iii),
as $\alpha_{r_p^B,p}$ is one of the leading Ulm-Kaplanski invariants of $B$.
This contradiction proves that 
$n\Z(p^{r_p^B})^{(\alpha_{{r_p^B},p})}=\{0\}$.
Therefore, $p^{r_p^B}$ must divide $n$. 
Since this holds 
for each 
$p\in P_B$,
this shows that
$n_B=\prod_{p\in P_B} p^{r_p^B}$ divides $n$. 
Since $n_B B=\{0\}$,
we conclude 
that $nB=\{0\}$.

\smallskip
(iii)$\to$(ii)
 Since $m_G G$ is \cd\ by Definition \ref{definition:of:m_G}, $|m_G G|=|m_G G|_d$.
Applying Lemma \ref{general:lemma:on:multiples:new}~(ii)
to the monotone cardinal function $\varphi(G)=|G|$,
 we get
$ |m_G G|_d=|G|_d$.
Furthermore,
$|G|_d=|H|$ by 
Lemma \ref{rem:cardinality:of:H}~(ii).
 Since $B\subseteq G$, we have $|m_G B|\le |m_G G|=|H|$.
Therefore, $n_B$ divides $m_G$ by 
(iv).
In particular, $n_B\le m_G$.
The reverse inequality $m_G\le n_B$ follows from Lemma \ref{m_G:is:a:lower:bound:for:n_B}.

\smallskip
(ii)$\to$(i) Suppose that $n_B=m_G$, yet
$n_B\neq b_G$.
Since $b_G\le n_B$ by Definition \ref{def:b_G}, 
we get
$b_G<n_{B}$. By \eqref{eq:33:fo},
there exists a Galindo-Macario decomposition $G=B'\oplus H'$ of $G$ such that $n_{B'}=b_G$. 
Since $m_G\le n_{B'}$ by Lemma \ref{m_G:is:a:lower:bound:for:n_B},
we get $m_G\le n_{B'}=b_G<n_B$, in contradiction with the assumption of (ii).
\end{proof}

\begin{corollary}
\label{algebraic:corollary}
Every abelian group $G$
has a (minimal) Galindo-Macario decomposition
$G=B\oplus H$ such that 
$b_G=n_B=m_G$.
In particular, the equality $b_G=m_G$ holds for every abelian group $G$.
\end{corollary}
\begin{proof}
Let $G=B\oplus H$ be a minimal Galindo-Macario decomposition which exists by Remark \ref{rem:minimal:decomposition}~(ii).
Then $b_G=n_B$ by Definition \ref{def:b_G},
and $n_B=m_G$ by the implication (i)$\to$(ii) of Theorem \ref{characterizations:of:minimality}.
\end{proof}

\section{Proofs of Theorems 
\ref{thm:6.2}
and \ref{no:GM}} 
\label{sec:second:group:of:proofs}
	
\begin{lemma}
\label{new:lemma}
Let $H$ be an
abelian group and $B$ be a bounded torsion abelian group
having 
the following property:
If $n\in\N$ and $|nB|\le |H|$, then $nB=\{0\}$.
If $B\oplus H$ admits a pseudocompact group topology,
then so does $B$.
\end{lemma}
\begin{proof}
Suppose that the direct sum $B\oplus H$ admits a pseudocompact group topology. Then, for every $n\in\N$, the group 
$
n(B\oplus H)=nB\oplus nH
$
admits a pseudocompact group topology by
Fact \ref{folklore}~(ii),
so its cardinality 
\begin{equation}
\label{eq:32:nm}
\tau_n=|n(B\oplus H)|=|nB\oplus nH|
\end{equation}
is either finite or an admissible cardinal
(Fact \ref{folklore}~(ii)).

Let $n\in\N$ be arbitrary. 
Consider two cases.

\smallskip
{\it Case 1.}
$|nB|>|H|$. 
Since $nH\subseteq H$, we have $|H|\ge |nH|$, so
$|nB|>|nH|$.
If the cardinality $|nB|$ of the group $nB$ is infinite, then 
$\tau_n=|nB|$ by \eqref{eq:32:nm} and $|nB|>|nH|$, so the cardinal $|nB|$ is admissible.

\smallskip
{\it Case 2.}
$|nB|\le |H|$. In this case $nB=\{0\}$ by
our assumption on $B$, so the cardinal $|nB|=1$ is finite in this case.

\smallskip
It follows from the above that, for every $n\in\N$, the cardinal $|nB|$ is either finite or admissible, 
so
$B$  
has
a pseudocompact group topology by Corollary
\ref{necessary:and:sufficient:condition:for:torsion:groups}.
\end{proof}

\noindent
{\bf Proof of Theorem  \ref{thm:6.2}:}
It follows from 
the implication (i)$\to$(iii) of 
Theorem \ref{characterizations:of:minimality}
that $B$ satisfies the assumption of  
Lemma \ref{new:lemma}.
Since $G=B\oplus H$ is pseudocompact, 
this lemma
implies that $B$ admits a pseudocompact 
group topology, and the
zero-dimensionality of this
topology 
follows from
Theorem \ref{psc:bounded:zero-dimensional}~(i).
\qed

\bigskip
\noindent
{\bf Proof of Theorem \ref{no:GM}:}
Let $G=B\oplus H$ be 
a minimal Galindo-Macario decomposition 
of $G$ (which exists by Remark \ref{rem:minimal:decomposition}~(ii)).
Note that $B=B[m]$ and 
\begin{equation}
\label{eq:mB=0}
mB=\{0\},
\end{equation}
 as 
$n_B$ is the order of $B$,
and $n_B=m$ by
the implication (i)$\to$(ii) of 
Theorem \ref{characterizations:of:minimality}.
Furthermore, $G[m]=B[m]\oplus H[m]=B\oplus H[m]$, so
\begin{equation}
\label{eq:29:zc}
nG[m]=nB\oplus n H[m]
\text{ for every }
n\in\N.
\end{equation}

\smallskip
The implication (iii)$\to$(ii) is trivial.

\smallskip
(ii)$\to$(i) 
Suppose that $n\in\N$ and 
$|nB|\le |H[m]|$. Since $H[m]\subseteq H$, we get 
$|H[m]|\le |H|$. This shows that $|nB|\le |H|$, so $nB=\{0\}$ by
the implication (i)$\to$(iii) of 
Theorem \ref{characterizations:of:minimality}.
This means that $B$ and $H[m]$ (taken is $H$ in Lemma \ref{new:lemma}) satisfy the assumption of Lemma \ref{new:lemma}.
Since $G[m]=B\oplus H[m]$ admits a pseudocompact group topology by (ii), from Lemma
\ref{new:lemma} we conclude that $B$ can be equipped with a pseudocompact group topology.

By (ii), \eqref{eq:mB=0} and Lemma \ref{powers:of:GM:decompositions}~(i), the quotient group $G/G[m]\cong mG=mH$ admits a pseudocompact group topology.  
Since $H$ is \cd, $H$ admits a 
pseudocompact group topology
by Theorem \ref{clever:lemma}.

Since both $B$ and $H$ admit a pseudocompact group topology, so does $G$ by Corollary \ref{equivalent:conditions:for:psc:topologizations}.

\smallskip
(i)$\to$(iii) Assume that $G$ admits a pseudocompact group topology. 
Then both $B$ and $H$ admit a pseudocompact group topology by 
Theorems \ref{thm:6.2}
and 
\ref{weight:divisible:part}.
It follows from Fact \ref{folklore}~(ii) that 
$mH$ admits a pseudocompact group topology.
If the \cd\ group
$mH$ is non-trivial,  
$mH$ admits a connected, locally connected pseudocompact group topology by Corollary 
\ref{structure:of:card-divisible:groups}.
If $mH=\{0\}$, then it is connected, locally connected and compact (so pseudocompact by \eqref{compact:implies:pseudocompact}) in its unique (group) topology. 
Since 
$mG=mH$ by \eqref{eq:mB=0} and Lemma \ref{powers:of:GM:decompositions}~(i),
$G/G[m]\cong mG$ admits a connected, locally connected pseudocompact group topology.

It remain only to show that the torsion group $G[m]$ admits a (zero-dimensional) pseudocompact group topology. To show this,
it suffices to check
that $G[m]$ satisfies two conditions of 
Corollary \ref{necessary:and:sufficient:condition:for:torsion:groups}.
Clearly, $G[m]$ is bounded.
Let $n\in\N$ be arbitrary. 
Consider two cases.

\smallskip
{\it Case 1.}
$|nB|>|H|$. 
Since $nH[m]\subseteq H$, we have $|H|\ge |nH[m]|$, so
$|nB|>|nH[m]|$.
If the cardinality $|nB|$ of the group $nB$ is infinite, then $|n G[m]|=|nB|$ by \eqref{eq:29:zc}. Since $B$ admits a pseudocompact group topology,
the cardinal $|nB|$ is admissible by Fact \ref{folklore}~(ii).
Hence,
the cardinal $|n G[m]|$ is admissible in this case.

\smallskip
{\it Case 2.}
$|nB|\le |H|$. In this case $nB=\{0\}$ 
 by 
the implication (i)$\to$(iii) of 
Theorem \ref{characterizations:of:minimality},
so
$nG[m]=n H[m]$ by \eqref{eq:29:zc}.
Since $m=n_B$,
and $n_B$ 
divides $n$ by 
the implication (i)$\to$(iii) of 
Theorem \ref{characterizations:of:minimality},
$m$ divides $n$.
Since $mH[m]=\{0\}$,
we get $n H[m]=\{0\}$,
so
the cardinal $|n G[m]|=|n H[m]|=1$ is finite in this case.

\smallskip
It follows from the above that, for every $n\in\N$, the cardinal $|n G[m]|$ is either finite or admissible. Therefore,
$G[m]$ admits a pseudocompact group topology by Corollary 
\ref{necessary:and:sufficient:condition:for:torsion:groups}.
This topology is zero-dimensional by Theorem \ref{psc:bounded:zero-dimensional}~(i).
\qed

\bigskip
\noindent
{\bf Statement on the use of AI:\/} The authors certify that this maniscript is completely ``AI-free''; that is, AI was {\em not\/} used in any capacity or form at any stages of its preparation, including (but not limited to) generation of main ideas, statements of results and their proofs, as well as proofreading.

\end{document}